\documentclass{amsart}
\copyrightinfo{2006}{American Mathematical Society}
\usepackage{mathrsfs,amscd}
\usepackage{amssymb}
\usepackage{amssymb, amsmath, mathrsfs}
\usepackage{amscd}
\usepackage{latexsym}
\usepackage{amssymb}
\usepackage{amsmath}
\usepackage{amsthm}
\usepackage{indentfirst}

\newtheorem{theorem}{Theorem}[section]

\theoremstyle{definition}
\newtheorem{definition}[theorem]{Definition}

\theoremstyle{remark}
\newtheorem{remark}[theorem]{Remark}

\numberwithin{equation}{section}
\allowdisplaybreaks
\begin{document}

\title{Explicit self-similiar singularity of Born-Infeld equation, membrane equation and space-like surfaces with vanishing mean curvature equation}

\author[W. Yan]{Weiping Yan}
\address{School of Mathematics, Xiamen University, Xiamen 361000, People's Republic of China}
\address{Laboratoire Jacques-Louis Lions, Université Pierre et Marie Curie (Sorbonne Université), 4, Place Jussieu, 75252 Paris, France.}\email{yanwp@xmu.edu.cn}
\thanks{The author is supported by NSFC No 11771359.}

\date{Dec 2017}

\begin{abstract}
In this short paper, we study the self-similar singularity phenomenon of zero mean curvature equation including Born-Infeld equation, membrane equation and space-like surfaces with vanishing mean curvature equation, which arises in string theory and geometric minimal surfaces theory. We show that there exist new explicit self-similar solutions for Born-Infeld equation, membrane equation and space-like surfaces with vanishing mean curvature equation in Minkowski space $R^{1+3}$ with respect to the radially symmetric case. 
\end{abstract}
\maketitle


\section{Introduction and main results}

The study of singularity is one of most important topics in physics and mathematics theory, which corresponds to a physical event, such as the solution (e.g. a physical flow field) changing topology, or the emergence of a new structure, such as a tip, cusp. It can also imply that some essential physics is missing from the equation in question, which should thus be supplemented with additional terms. There has been discovered that the behavior of string theory in spacetimes that develop singularities \cite{W}. In string theory, relativistic strings and membranes correspond to the timelike minimal surface which can be described by Born-Infeld equation and membrane equation. 
The Born-Infeld equation was first established by Born and Infeld \cite{Born0,Born} to describe nonlinear electrodynamics, a generalization of the linear Maxwell equations. 
It is also appears in a geometric nonlinear theory of electromagnetism, and seen as the equation of graphs with zero mean curvature over a domain of the timelike $tx$-plane in Lorentz-Minkowski $\mathbb{L}^3(t,x,y)$. Membrane equation is a higher dimensional version of Born-Infeld equation.
One can see \cite{Che} for a detailed introduction to this field.
Generic singularity of them has attracted many interesting works \cite{D1,D2,hop,Hop2,V}. 
Since there is much difficulty to find an explicit solution of them, only approximation of solutions of it has been obtained in \cite{Fe,Fe1,Fe3}. To author's knowledgement, there is few result on explicit singularity of this kind of equations. Luckily, we find some explicit singularity (self-similar singularity) of them (Born-Infeld equation, space-like surfaces with vanishing mean curvature equation and membrane equation) through observation of some structure of them in this paper.

We consider the Born-Infeld equation, membrane equation and space-like surfaces with vanishing mean curvature equation as follows. 

The classical Born-Infeld equation:
\begin{equation}\label{E1-1}
u_{tt}(1+u_x^2)-u_{xx}(1-u_t^2)=2u_tu_xu_{tx}, \quad (t,x)\in\mathbb{R}^+\times\mathbb{R}.
\end{equation}

The radially symmetric membrane equation ($r=|x|$, $x\in\mathbb{R}^3$)
\begin{equation}\label{E1-3}
u_{tt}-u_{rr}-\frac{u_r}{r}+u_{tt}u_r^2+u_{rr}u_t^2-2u_tu_ru_{tr}+\frac{1}{r}u_ru_t^2-\frac{1}{r}u_r^3=0,
 \quad (t,r)\in\mathbb{R}^+\times\mathbb{R}^+.
\end{equation}

The space-like surfaces with vanishing mean curvature equation
\begin{equation}\label{E1-4}
u_{xx}(1-u_y^2)+u_{yy}(1-u_x^2)+2u_xu_yu_{xy}=0,\quad (x,y)\in\mathbb{R}\times\mathbb{R}.
\end{equation}

Here all $u$ in above equations are \textbf{scalar} unknown functions.

The classical Born-Infeld equation and membrane equation can be derived from the Lagrange action. In general \cite{Hop1},
let $\mathcal{M}$ be a timelike $(M+1)$-dimensional hypersurface, and $(\mathbb{R}^{D},g)$ be a $D$-dimensional Minkowski space, and $g$ be the Minkowski metric with $g(\partial_t,\partial_t)=1$. At any time $t$, the spacetime volume in $\mathbb{R}^{D}$ of timelike hypersurface $\mathcal{M}$ can be described as a \textbf{graph} over $\mathbb{R}^M$, which satisfies
\begin{equation}\label{E1-0}
\mathcal{S}(u)=\int_{\mathbb{R}^M}\sqrt{1-|\partial_t u|^2+|\nabla u|^2}d^Mxdt.
\end{equation}
Critical points of action integral (\ref{E1-0}) give rise to submanifolds $\mathcal{M}\subset\mathbb{R}^{D}$ with
vanishing mean curvature, i.e. time-like extremal hypersurfaces (membrane). The Euler-Lagrange equation of (\ref{E1-0}) takes the form
\begin{equation}\label{ENNN1-1}
\partial_t\left(\frac{\partial_t u}{\sqrt{1-|\partial_t u|^2+|\nabla_xu|^2}}\right)-\sum_{i=1}^M\partial_{x_i}\left(\frac{\partial_{x_i}u}{\sqrt{1-|\partial_t u|^2+|\nabla_xu|^2}}\right)=0.
\end{equation}

The Born-Infeld equation and radially symmetric membrane equation are \textbf{graph equations} (i.e. for Born-Infeld equation, the graph is $(t,x,u(t,x))\in\mathbb{R}^3$) of (\ref{E1-0}).
When $M=1$ and $D=3$, equation (\ref{ENNN1-1}) is equivalent to the classical Born-Infeld equation (\ref{E1-1}). Here $u:\mathbb{R}\times\mathbb{R}\rightarrow\mathbb{R}$.

When $M=2$ and $D=4$, equation (\ref{ENNN1-1}) is equivalent to membrane equation (also called the timellike minimal surface equation)
\begin{equation}\label{E1-2}
(1-u_{\alpha}u^{\alpha})\square u+u^{\beta}u^{\alpha}u_{\alpha\beta}=0,
\end{equation}
where $\forall\alpha,\beta=0,1,2,\ldots,M$, $u_{\alpha}=\frac{\partial u}{\partial x^{\alpha}}$, $z_{\alpha\beta}=\frac{\partial^2 u}{\partial x^{\alpha}\partial x^{\beta}}$ and $\square u=u_{\alpha\beta}g^{\alpha\beta}$.

Especially, let $r=|x|$, then equation (\ref{E1-2}) is reduced into the radially symmetric membrane equation
\begin{equation*}
u_{tt}-u_{rr}-\frac{u_r}{r}+u_{tt}u_r^2+u_{rr}u_t^2-2u_tu_ru_{tr}+\frac{1}{r}u_ru_t^2-\frac{1}{r}u_r^3=0.
\end{equation*}

Different with above Lagrange action (\ref{E1-0}) derivations, another interesting zero-mean curvature equation is space-like surfaces with vanishing mean curvature equation (\ref{E1-4}), which is graphs with zero mean curvature over a domain of the spacelike $xy$-plane in Lorentz-Minkowski $\mathbb{L}^3(t,x,y)$.

It is easy to see that Born-Infeld equation (\ref{E1-1}) and the symmetric membrane equation (\ref{E1-3})  exhibit the following scaling invariance for any $\lambda>0$,
\begin{equation}\label{E1-5}
u(t,x)\mapsto u_{\lambda}(t,x)=\lambda^{-1} u(\lambda t,\lambda x).
\end{equation}
Under this scaling the conserved energy of equation (\ref{E1-1}) and (\ref{E1-3}) 
\begin{equation*}
E_1(u)=\int_{0}^{\infty}(\frac{1}{2}u_t^2+\frac{1}{2}u_x^2+F_1(u_t,u_x))xdx,
\end{equation*}
and
\begin{equation*}
E_2(u)=\int_{0}^{\infty}(\frac{1}{2}u_t^2+\frac{1}{2}u_r^2+F_2(u_t,u_r))rdr,
\end{equation*}
where 
\begin{eqnarray*}
&&F'_1(u_t,u_x)=u_{tt}u_x^2+u_{xx}u_t^2-2u_tu_xu_{tx},\\
&&F'_2(u_t,u_r)=u_{tt}u_r^2+u_{rr}u_t^2-2u_tu_ru_{tr}+\frac{1}{r}u_ru_t^2-\frac{1}{r}u_r^3.
\end{eqnarray*}
can be transformed as
\begin{equation*}
E_1(u_{\lambda})=\lambda E_1(u),~~E_2(u_{\lambda})=\lambda E_2(u).
\end{equation*}
Moverover, it is a mass conservation dynamics, it follows from (\ref{ENNN1-1}) that the quantity
\begin{equation*}
\int\left(\frac{\partial_t u}{\sqrt{1-|\partial_t u|^2+|\nabla_xu|^2}}\right)~is~conserved~along~the~dynamics.
\end{equation*}
Since equations (\ref{E1-1}) and (\ref{E1-3}) are energy supercritical, so one expects smooth finite energy initial data to lead to finite time blow up, and the blow up rate is like the self-similar blow up solution. 
Eggers and Hoppes \cite{Hop1} gave a detail discussion on the existence of self-similar blow up solutions (not explicit self-similar solutions) to the Born-Infeld equation (\ref{E1-1}). 
They showed that above equation has self-similar solutions
\begin{equation*}
u(t,x)=u_0-\hat{t}+\hat{t}^{a}h(\frac{x}{\hat{t}^b})+\ldots,
\end{equation*}
where $\hat{t}=t_0-t$ and $h(x)\varpropto A_{\pm}x^{\frac{2a}{a+1}}$ for $x\rightarrow\pm\infty$. To the higher dimension case, they showed that the radially symmetric membranes equation (\ref{E1-2}) has a self-similar solutions
\begin{equation*}
u(t,x)=-\hat{t}+\hat{t}^{a}h(\frac{x-x_0}{\hat{t}^b})+\ldots,
\end{equation*}
by analyzing the eikonal equation
\begin{equation*}
1-u_t^2+u_x^2=0.
\end{equation*}
Meanwhile, the swallowtail singularity was also been given by parametric the string solution in \cite{Egg}. One can see \cite{Che0,Fe2,Fe,Fe1,Fe3,hop,Egg,Hop2,Lin,It,tian} for more disscusion on the existence of solutions for the Born-Infeld equation and membrane equation. Recently, Alejo and Mu$\tilde{n}$oz \cite{Alejo} obtained a sharp nonlinear scattering result for Born-Infeld equation (\ref{E1-2}).

In this paper, we find that there are explicit self-similar blow up solutions to Born-Infeld equation (\ref{E1-1}), the symmetric membrane equation (\ref{E1-3}) and space-like surfaces with vanishing mean curvature equation (\ref{E1-4}). Moreover, the mode unstable of explicit self-similar solution of equation (\ref{E1-3}) is given. To author's knowlege, this is the only result on explicit self-similar singularity solutions for equations (\ref{E1-1}), (\ref{E1-3})-(\ref{E1-4}).

More precisely, the main result is the following.
\begin{theorem}
The Born-Infeld equation (\ref{E1-1}), the radial symmetric membrane equation (\ref{E1-3}) and space-like surfaces with vanishing mean curvature equation (\ref{E1-4}) have explicit self-similar blow up solutions
$$
u_k(t,x)=k\ln(\frac{T-t+x}{T-t-x}),\quad |x|\leq T-t,\quad t\in[0,T),
$$
$$
u_{\pm}(t,r)=\pm(T-t)\sqrt{1-(\frac{r}{T-t})^2},\quad r\in[0,T-t],\quad t\in[0,T),
$$
and
$$
u_k(x,y)=k\ln|\frac{y}{T-x}+(1+\frac{y^2}{(T-x)^2})^{\frac{1}{2}}|,
$$
respectively, where $T$ is a positive constant and $k$ is an arbitrary constant in $\mathbb{R}\setminus \{0\}$.

Moreover, there are unstable mode $\nu=4$ and stable mode $\nu=-1$ for the linearized symmetric membrane equation (\ref{E1-3}) around $u_{\pm}(t,r)$.
\end{theorem}

\begin{remark}
From the form of explicit self-similar solution of Born-Infeld equation in above theorem, one can see that it is defined in the lightcone
$$
\mathcal{L}_T:=\{(t,x) | |x|<T-t,\quad t\in[0,T)\}.
$$
Hence $x$ and $r$ are finite range which is determined by the form of explicit solutions, and those explicit self-similar solutions of Born-Infeld equation and membrane equation describe the local behavior of singularity in the lightcone. 

It is easy to see that the shape of explicit self-similar solution on the radial symmetric membrane equation is a sphere in $\mathbb{R}^{1,3}$.
\end{remark}

\section{Proof of Theorem 1.1}
\subsection{Born-Infeld equation (\ref{E1-1})}
Let parameter $T$ be a positive constant.
Since self-similar solutions are invariant under the scaling (\ref{E1-5}),
we introduce the similarity coordinates
\begin{equation*}
\tau=-\log(T-t),~~~~\rho=\frac{x}{T-t},
\end{equation*}
then we denote by
\begin{equation*}
u(t,x)=v(-\log(T-t),\frac{x}{T-t}),
\end{equation*}
and noticing 
\begin{eqnarray*}
&&u_t(t,x)=e^{\tau}(v_{\tau}+\rho v_{\rho}),\\
&&u_{tt}(t,x)=e^{2\tau}(v_{\tau\tau}+v_{\tau}+2\rho v_{\rho}+2\rho v_{\tau\rho}+\rho^2v_{\rho\rho}),\\
&&u_x(t,x)=e^{\tau}v_{\rho},\\
&&u_{xx}(t,x)=e^{2\tau}v_{\rho\rho},\\
&&u_{tx}(t,x)=e^{2\tau}(v_{\tau\rho}+v_{\rho}+\rho v_{\rho\rho})
\end{eqnarray*}
equation (\ref{E1-1}) is transformed into an one dimensional quasilinear wave equation
\begin{eqnarray}\label{YN1-1}
v_{\tau\tau}-(1-\rho^2)v_{\rho\rho}&+&v_{\tau}+2\rho v_{\rho}+2\rho v_{\tau\rho}+e^{2\tau}v_{\rho}^2(v_{\tau\tau}+v_{\tau}+2\rho v_{\rho}+2\rho v_{\tau\rho}+\rho^2v_{\rho\rho})\nonumber\\
&&+e^{2\tau}(v_{\tau}+\rho v_{\rho})^2v_{\rho\rho}-2e^{2\tau}v_{\rho}(v_{\tau}+\rho v_{\rho})(v_{\rho}+\rho v_{\rho\rho}+v_{\tau\rho})=0.
\end{eqnarray}
The steady equation of wave equation (\ref{YN1-1}) is 
\begin{equation*}
(\rho^2-1)v_{\rho\rho}+2\rho v_{\rho}=0,
\end{equation*}
which is an ODE. Direct computation shows that it has a family of solutions
\begin{equation*}
v(\rho)=k\ln\frac{1+\rho}{1-\rho},
\end{equation*}
where $k$ is an arbitrary constant in $\mathbb{R}\setminus \{0\}$.

Obviously, the domain of $\rho$ is 
\begin{equation*}
\{\rho| |\rho|<1\}.
\end{equation*}
Hence Born-Infeld equation (\ref{E1-1}) has a family of explicit self-similar solutions
\begin{equation*}
u_k(t,x)=k\ln(\frac{T-t+x}{T-t-x}),~~k\in \mathbb{R}\setminus \{0\}.
\end{equation*}
It is easy to see that
\begin{equation*}
\partial_xu_k|_{x=0}=\frac{k}{T-t}\rightarrow+\infty,~~as~~t\rightarrow T^{-}.
\end{equation*}

\subsection{The symmetric  membrane equation (\ref{E1-3})}
Self-similar solutions are invariant under the scaling (\ref{E1-5}), so let
\begin{eqnarray*}
&&\rho=\frac{r}{T-t},\\
&&u(t,r)=(T-t)\phi(\rho),
\end{eqnarray*}
where $T$ is a positive constant, which is a parameter.

Inseting this ansatz into equation (\ref{E1-3}) by noticing
\begin{eqnarray*}
&&\partial_tu(t,r)=-\phi(\rho)+\rho\phi'(\rho),\\
&&\partial_{tt}u(t,r)=(T-t)^{-1}\rho^2\phi''(\rho),\\
&&\partial_ru(t,r)=\phi'(\rho),\\
&&\partial_{rr}u(t,r)=(T-t)^{-1}\phi''(\rho),\\
&&\partial_{tr}u(t,r)=(T-t)^{-1}\rho\phi''(\rho),
\end{eqnarray*}
one gets a quasilinear ordinary differential equation
\begin{equation}\label{E2-1}
\rho(1-\rho^2)\phi''+\phi'-\phi'\phi^2+2\rho\phi(\phi')^2-\rho\phi''\phi^2+(1-\rho^2)(\phi')^3=0.
\end{equation}
We are interested in smooth solutions in the backward light-cone of the blow up point $(t,r)=(T,0)$, i.e. $\rho$ in the closed interval $[0,1]$. This solution is such that
\begin{equation}\label{E2-1R1}
\frac{\partial^n}{\partial r^n}\phi(\frac{r}{T-t})|_{r=0}=(T-t)^{-n}\frac{d^n\phi}{dy^n}(0),
\end{equation}
where the $n$-th derivative is even, there is $\frac{d^n\phi}{dy^n}(0)\neq0$, this means it is diverges as $t\rightarrow T$. When $n$ is odd, there is $\frac{d^n\phi}{dy^n}(0)=0$. This condition is different with wave map (e.g. see \cite{BB}).
One can see each self-similar solution $\phi(\rho)\in\mathbb{C}^{\infty}[0,1]$ describes a singularity developing in finite time from smooth initial data. Since (\ref{E2-1}) is quasilinear ODE, there are some difficulties to prove the existence of smooth of solutions.
But from the structure of nonlinear term in equation (\ref{E2-1}), it can be rewritten as
\begin{equation}\label{E2-1R}
\rho(1-\rho^2-\phi^2)\phi''+\phi'-\phi'\phi^2+2\rho\phi(\phi')^2+(1-\rho^2)(\phi')^3=0.
\end{equation}
We found that if we let
\begin{equation*}
1-\rho^2-\phi^2=0,
\end{equation*}
then there holds
\begin{equation*}
\phi'-\phi'\phi^2+2\rho\phi(\phi')^2+(1-\rho^2)(\phi')^3=0.
\end{equation*}
Hence we obtain two explicit solutions
\begin{equation}\label{E2-2}
\phi(\rho)=\pm\sqrt{1-\rho^2},
\end{equation}
which satisfies condition (\ref{E2-1R1}).

Consequently, two explicit self-similar solutions of equation (\ref{E1-3}) are
\begin{equation}\label{E2-3}
u_{\pm}(t,r)=\pm(T-t)\sqrt{1-(\frac{r}{T-t})^2},
\end{equation}
which exhibit smooth for all $0<t<T$, but which break down at $t=T$ in the sense that
\begin{eqnarray*}
\partial_{rr}u_{\pm}(t,r)|_{r=0}&=&\pm((T-t)^2-r^2)^{-\frac{1}{2}}|_{r=0}\pm r^2((T-t)^2-r^2)^{-\frac{3}{2}}|_{r=0}\nonumber\\
&=&\pm\frac{1}{T-t}\rightarrow+\infty,~~as~~t\rightarrow T^{-},
\end{eqnarray*}
and the dynamical behavior of them are as attractors.

On the other hand, from the form of $u_{\pm}(t,r)$ in (\ref{E2-3}), it requires that
\begin{equation*}
1-(\frac{r}{T-t})^2\geq0.
\end{equation*}
So we consider the dynamical behavior of self-similar solutions in
of the backward lightcone
\begin{equation*}
\mathcal{B}_T:=\{(t,r):t\in(0,T),~~r\in[0,T-t]\}.
\end{equation*}

\begin{remark}
Obviouly, self-similar solutions (\ref{E2-3}) are cycloids. The sphere begins to expand until it starts to shrink and eventually collapses to a point in a finite time $\tilde{T}=T-r_0$, i.e. $u_{\pm}(\tilde{T},r_0)=0$.
Here $r_0$ is a fixed positive constant in the backward lightcone.
\end{remark}

We introduce the similarity coordinates
\begin{equation}\label{E2-5}
\tau=-\log(T-t),~~~~\rho=\frac{r}{T-t},
\end{equation}
then we denote by
\begin{equation*}
\tilde{v}(\tau,\rho)=e^{\tau}u(T-e^{-\tau},\rho e^{-\tau}),
\end{equation*}
equation (\ref{E1-3}) is transformed into
\begin{eqnarray}\label{E2-4}
\tilde{v}_{\tau\tau}-\tilde{v}_{\tau}&-&(1-\rho^2)\tilde{v}_{\rho\rho}-\frac{1}{\rho}\tilde{v}_{\rho}+2\rho \tilde{v}_{\tau\rho}+\tilde{v}_{\rho}^2(\tilde{v}_{\tau\tau}+\tilde{v}_{\tau}-2\tilde{v})+\tilde{v}_{\rho\rho}
(\tilde{v}-\tilde{v}_{\tau})^2\nonumber\\
&&-2\tilde{v}_{\rho}\tilde{v}_{\tau\rho}(\tilde{v}_{\tau}-\tilde{v})+\frac{1}{\rho}\tilde{v}_{\rho}(\tilde{v}_{\tau}-\tilde{v})^2+\frac{1}{\rho}(\rho^2-1)\tilde{v}_{\rho}^3=0.
\end{eqnarray}
In the similarity coordinates (\ref{E2-5}), the blow up time $T$ is changed to $\infty$. So the stability of blow up solutions of the radially symmetric membranes equation  (\ref{E1-3}) as $t\rightarrow T^{-}$ is transformed into the asymptotic stability of quasilinear wave equation (\ref{E2-4}) as $\tau\rightarrow\infty$. Now we consider the problem of mode stable or unstable of linear equation of (\ref{E2-4}). In spirt of the method of wave map (e.g. see \cite{BB,Cos2}), set solution of (\ref{E2-4}) takes the form
\begin{eqnarray}\label{E2-4R}
\tilde{v}(\tau,\rho)&=&\phi(\rho)+v(\tau,\rho),
\end{eqnarray}
where $\phi(\rho)$ is defined in (\ref{E2-2}).

Inserting this ansatz (\ref{E2-4R}) into equation (\ref{E2-4}), $v(\tau,\rho)$ satisfies one dimensional equation
\begin{equation}\label{E2-6r}
\begin{array}{lll}
(1+\phi'^2)v_{\tau\tau}&-&(1-\phi'^2+2\phi''\phi+\frac{2}{\rho}\phi')v_{\tau}+2(\phi\phi'+\rho)v_{\tau\rho}-(1-\rho^2-\phi^2)v_{\rho\rho}\\
&&-\frac{1}{\rho}(1+4\rho\phi'\phi-3(\rho^2-1)\phi'^2-\phi^2)v_{\rho}+\frac{1}{\rho}(-2\rho\phi'^2+2\rho\phi''\phi+2\phi'\phi)v\\
&&-2\phi v_{\rho}^2+v_{\rho}(v_{\tau\tau}+v_{\tau}-2v)(v_{\rho}+2\phi')+\phi''(v-v_{\tau})^2\\
&&+v_{\rho\rho}[(v-v_{\tau})^2+2\phi(v-v_{\tau})]-2\phi'v_{\tau\rho}(-v+v_{\tau})+2\phi v_{\rho}v_{\rho\tau}\\
&&-2v_{\rho}v_{\rho\tau}(-v+v_{\tau})+\frac{1}{\rho}v_{\rho}[(v-v_{\tau})^2+2\phi(v-v_{\tau})]+\frac{1}{\rho}\phi'(v-v_{\tau})^2\\
&&+(\rho-\frac{1}{\rho})(v_{\rho}^3+3\phi'v_{\rho}^2)=0.
\end{array}
\end{equation}
From the exact form of $\phi(\rho)$ in (\ref{E2-2}), it has
\begin{eqnarray*}
&&1-\rho^2-\phi^2=0,\\
&&\phi\phi'+\rho=0.
\end{eqnarray*}
Thus one dimensional equation (\ref{E2-6r}) is loss of hyperbolicity, which takes the form
\begin{equation}\label{E2-7r}
\begin{array}{lll}
v_{\tau\tau}&+&3v_{\tau}-4v-2(1-\rho^2)^{\frac{1}{2}}v_{\rho}^2+v_{\rho}(v_{\tau\tau}+v_{\tau}-2v)(v_{\rho}-2\rho(1-\rho^2)^{-\frac{1}{2}})\\
&&-(1-\rho^2)^{-\frac{3}{2}}(v-v_{\tau})^2+v_{\rho\rho}[(v-v_{\tau})^2+2(1-\rho^2)^{\frac{1}{2}}(v-v_{\tau})]\\
&&+2\rho(1-\rho^2)^{-\frac{1}{2}}v_{\tau\rho}(-v+v_{\tau})+2(1-\rho^2)^{\frac{1}{2}} v_{\rho}v_{\rho\tau}-2v_{\rho}v_{\rho\tau}(-v+v_{\tau})\\
&&+\frac{1}{\rho}v_{\rho}[(v-v_{\tau})^2+2(1-\rho^2)^{\frac{1}{2}}(v-v_{\tau})]-(1-\rho^2)^{-\frac{1}{2}}(v-v_{\tau})^2\\
&&+(\rho-\frac{1}{\rho})(v_{\rho}^3-3\rho(1-\rho^2)^{-\frac{1}{2}}v_{\rho}^2)=0.
\end{array}
\end{equation}
then we linearize (\ref{E2-6r}) at $v(\tau,\rho)=e^{\nu\tau}u_{\nu}$ leading to a eigenvalue problem
\begin{equation}\label{E2-4R1}
(\nu^2+3\nu-4)u_{\nu}=0.
\end{equation}
As in \cite{Cos2, Cos3}, we introduce defnitions of mode stable and unstable of solution $u_{\nu}$ to equation (\ref{E2-4R1}).

\begin{definition}
A non-zero smooth solution $u_{\nu}$ of (\ref{E2-4R1}) is called mode stable if $Re~\nu<0$ holds. The eigenvalue $\nu$ is called a stable eigenvalue. Otherwise, if $Re~\nu\geq0$ is called mode unstable of non-zero smooth solution $u_{\nu}$ of (\ref{E2-4R1}). $\nu$ is called an unstable eigenvalue.
\end{definition}

From (\ref{E2-4R1}), it is easy to check that $\nu=4$ and $\nu=-1$ are two eigenvalues of linear operator to $(\ref{E2-7r})$.
By the definition of mode unstable, we know that two explicit self-similar solutions $u_{\pm}(t,r)$ in (\ref{E2-3}) of time-like extremal hypersurfaces equation (\ref{E1-2}) are mode unstable.

\subsection{Space-like surfaces with vanishing mean curvature equation (\ref{E1-4})}
Introduce the similarity coordinates
\begin{equation*}
\tau=-\log(T-x),~~~~\rho=\frac{y}{T-x},
\end{equation*}
then we denote by
\begin{equation*}
u(x,y)=v(-\log(T-x),\frac{y}{T-x}),
\end{equation*}
and noticing 
\begin{eqnarray*}
&&u_x(x,y)=e^{\tau}(v_{\tau}+\rho v_{\rho}),\\
&&u_{xx}(x,y)=e^{2\tau}(v_{\tau\tau}+v_{\tau}+2\rho v_{\rho}+2\rho v_{\tau\rho}+\rho^2v_{\rho\rho}),\\
&&u_y(x,y)=e^{\tau}v_{\rho},\\
&&u_{yy}(x,y)=e^{2\tau}v_{\rho\rho},\\
&&u_{xy}(x,y)=(v_{\tau\rho}+v_{\rho}+\rho v_{\rho\rho})
\end{eqnarray*}
equation (\ref{E1-1}) is transformed into an one dimensional quasilinear elliptic equation
\begin{eqnarray}\label{YNN1-1}
v_{\tau\tau}+(1+\rho^2)v_{\rho\rho}&+&v_{\tau}+2\rho v_{\rho}+2\rho v_{\tau\rho}-e^{2\tau}v_{\rho}^2(v_{\tau\tau}+v_{\tau}+2\rho v_{\rho}+2\rho v_{\tau\rho}+\rho^2v_{\rho\rho})\nonumber\\
&&-e^{2\tau}(v_{\tau}+\rho v_{\rho})^2v_{\rho\rho}+2e^{2\tau}v_{\rho}(v_{\tau}+\rho v_{\rho})(v_{\rho}+\rho v_{\rho\rho}+v_{\tau\rho})=0.
\end{eqnarray}
The steady equation of elliptic equation (\ref{YNN1-1}) is 
\begin{equation*}
(\rho^2+1)v_{\rho\rho}+2\rho v_{\rho}=0,
\end{equation*}
which is an ODE. Direct computation shows that it has a family of solutions
\begin{equation*}
v(\rho)=k\ln|\rho+(1+\rho^2)^{\frac{1}{2}}|,
\end{equation*}
where $k$ is an arbitrary constant in $\mathbb{R}\setminus \{0\}$.

Space-like surfaces with vanishing mean curvature equation (\ref{E1-4}) has a family of explicit self-similar solutions
\begin{equation*}
u_k(x,y)=k\ln|\frac{y}{T-x}+(1+\frac{y^2}{(T-x)^2})^{\frac{1}{2}}|,~~k\in \mathbb{R}\setminus \{0\}.
\end{equation*}
It is easy to see that
\begin{equation*}
\partial_yu_k|_{y=0}=\frac{k}{T-x}\rightarrow+\infty,~~as~~x\rightarrow T^{-}.
\end{equation*}

\section{Discussion}
We have obtained some explicit self-similar singularities on $M$-dimensional timelike minimal surface in Minskowsi space $\mathbb{R}^{D}$ with $M=1,2$ and $D=2,3$. In fact, we think that there are many other singularities formation of Born-Infeld equation and membrane equation in Minkowski $\mathbb{R}^D$ with $D=2,3$. But they maybe have not an explicit mathematic formular. 
For instance, one can see the work of Eggers and Hoppe \cite{Egg}. At present, there is few paper concerning with the singluarity of  membrane equation in higher dimensional $(M\geq3,~D\geq4)$ Minkowski space except the work of \cite{Egg}. Hoppe \cite{Hop3} tells me that lightlike self-similar solutions obtained in Theorem 1.1. are also lightlike self-similar solutions of membrane equation in higher dimension $M\geq3,~D\geq4$.
Some meaningful problems are to consider the linear and nonlinear stability or instability of explicit self-similar solutions which are found in Theorem 1.1. In fact, we have shown there is a positve eigenvalue in the linearied symmetric membrane equation in Theorem 1.1. Recently, we have obtained there is a feedback control and stability of explicit self-similar solutions to the symmetric membrane equation (\ref{E1-3}) in \cite{Yan}. Consequently, it gives a way to stabilize two explicit self-similar solutions of equation (\ref{E1-3}). Furthermore, for the minimal surface equation
\begin{equation*}
u_{xx}(1+u_y^2)+u_{yy}(1+u_x^2)-2u_xu_yu_{xy}=0,
\end{equation*}
by Calabi's correspondence, it has a family of explicit solutions
\begin{equation*}
u_k(x,y)=k\ln|\frac{y}{T-x}+(1+\frac{y^2}{(T-x)^2})^{\frac{1}{2}}|,
\end{equation*}
which is the same with explicit solutions of space-like surfaces with vanishing mean curvature equation (\ref{E1-4}).
This makes us giving new gravitational instantons \cite{Yan3}.

Here we list some problems which will be considered in the furture.

i) Can we prove nonlinear stability or instability of explicit self-similar solutions which are obtained in Theorem 1.1 for Born-Infeld equation (\ref{E1-1}), space-like surfaces with vanishing mean curvature equation (\ref{E1-4})?

ii) Can we prove nonlinear instability of explicit self-similar solutions which are obtained in Theorem 1.1 for equation (\ref{E1-3})?

\textbf{Acknowledgments.} The author express his thanks to Prof. Alejo for sending his interesting paper \cite{Alejo} to me.

\end{document}